\documentclass[oneside]{article}
\usepackage{amsmath, graphicx, enumerate, xurl, tabularx, multirow, hyperref}
\usepackage[utf8]{inputenc}
\usepackage{mathptmx}
\usepackage{amssymb} 
\usepackage{devanagari} 
\title{Sa\.ngamagr\=ama M\=adhava: An Updated Biography}
\author{{\bf V. N. Krishnachandran}\\ Principal (retired), Govt. Victoria College, Palakkad, Kerala\\
(Formerly, Professor of Mathemaics, Govt. Engg. College, Thrissur)\\
Email: {\tt krishnachandranvn@gmail.com}}
\date{}
\begin{document}
\maketitle
\begin{abstract}
This paper presents an updated biography of Sa\.ngamagr\=ama M\=adhava incorporating viewpoints expounded by scholars in the recent past and collecting together in one place more details about his works culled from recent researches into his contributions to mathematics and astronomy. One major updation  is with regard to the geographical location in Kerala where M\=adhava flourished. Other updations include observations on M\=adhava's algorithm to compute the numerical values of the sine and cosine functions and on his value of the mathematical constant $\pi$. 
\end{abstract}
\tableofcontents
\newpage
%
\begin{center}
\includegraphics[height=7cm]{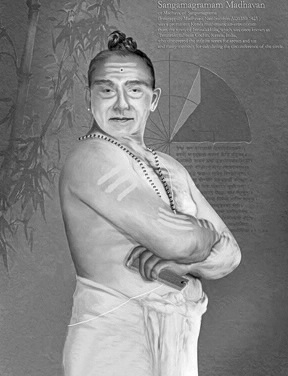}\\[2mm]
Sa\.ngamagr\=ama M\=adhava (c. 1340--1425)\footnote{This digital image of M\=adhava  is only  an artist’s impressions of what he might have looked like. It was drawn up with inputs provided by people, who some believe are M\=adhava's descendants, and released in 2014 by Madhava Ganitha Kendram, a Kochi based voluntary association working to revive his works.}
\end{center}
\section{Introduction}
Even though Sa\.ngamagr\=ama M\=adhava is a legendary figure in the mathematical and astronomical traditions of Kerala, practically very little is known with absolute certainty about \index{Sa\.ngamagr\=ama M\=adhava}\index{M\=adhava} the person M\=adhava, his life and his times. 
In these traditions he is more renowned as an astronomer and as an astrologer than as a mathematician. For example, in the encyclopedic five-volume history of Kerala literature ({\em Kera\d la S\=ahitya Caritra\d m})\index{Kera\d la S\=ahitya Caritra\d m@{\em Kera\d la S\=ahitya Caritra\d m}} written  by Ulloor S. Parameswara Aiyer\index{Ulloor S. Parameswara Aiyer} and published during 1953--1957, the author has referred to M\=adhava as a {\em Daivaj\~na}\index{Daivaj\~na} (an astrologer/astronomer) (see \cite{Ulloor} p. 98). 
Incidentally, M\=adhava's name appeared in print in a language other than Sanskrit for the first time only in 1948  in a Malayalam commentary of {\em Yuktibh\=a\d s\=a} by Ramavarma (Maru) Tampuran\index{Ramavarma (Maru) Tampuran} and
Akhilesvara Ayyar\index{Akhilesvara Ayyar} (see \cite{YuktiSayahna} p. xii).\footnote{M\=adhava's name appears in the Sanskrit astronomical literature of Kerala some of which were printed and published much earlier than this work. See, for example, \cite{Aryabhatiya1930}.} The much admired and much talked about C. M. Whish's paper on Kerala mathematics published in 1834 had no reference to M\=adhava (see \cite{Whish1834}). 
It was only during the second half of the twentieth century M\=adhava began to be recognised and adored as one of the greatest mathematicians Indian mathematical tradition had produced and the importance of his contributions to mathematics in the history of mathematics across the world began to be appreciated. 

\section{M\=adhava: The man and the place where he lived}
There are several references in the works of the astronomers and mathematicians of the Kerala school indicating that M\=adhava hailed from Sa\.ngama\-gr\=ama (see \cite{Sarma1958}). According to Ulloor S. Paramesvara Aiyer, the great historian of Kerala literature, there are some local legends to the effect that Sa\.ngamagr\=ama M\=adhava was ``Iri\~n\~n\=a\d tappa\d l\d li M\=adhavan Namp\=utiri'' hailing from  Sa\.ngamagr\=ama which is the town of Iri\~n\~n\=alakku\d da\index{Iri\~n\~n\=alakku\d da}, about 60 km north of Cochin and some 70 km south of the Ni\d la river (see \cite{Ulloor} p.98, \cite{Divakaran2008} p. 286,  \cite{Sarma1972} p. 51).\footnote{Ulloor S. Paramesvara Aiyer has also mentioned another legend to the effect that Sa\.ngamagr\=ama M\=adhava was a {\em V\=ariyar} belonging to the {\em Thekke\d tattu V\=ariya\d m} family also in Iri\~n\~n\=alakku\d da. However Aiyer, without citing any reasons whatsoever, discounted this legend as totally unreliable.} Further, according to Aiyar, these legends are reliable. 
However, as we shall see below, recent scholarship has contested this commonly accepted view.  Other than these legends and some highly stretched interpretations of place names and the connection between the house name of M\=adhava and the house name of a local Namp\=utiri\footnote{A Namp\=utiri is a Malayalam speaking Brahmin native to what is now the state of Kerala, India, where they constituted part of the traditional feudal elite. They were distinguished from other Brahmin sects by rare practices such as the adherence to \'Srauta ritualism, the P\=urva-M\=\i m\=a\d ms\=a school of Hindu philosophy and orthodox traditions, as well as many idiosyncratic customs that are unique among Brahmins, including primogeniture.}\index{Namp\=utiri} 
 family there appears to be no other concrete evidence to connect M\=adhava with a Namp\=utiri family in Irinj\=alakku\d da. 
 
In this section, we first examine an important Malayalam document that sheds some light on M\=adhava the man and then scrutinise in dome detail some of the terms in the document. We conclude the section with a  discussion on what modern scholarship has to say about the place name ``Sa\d ngamagr\=ama''. 

\subsection{An important Malayalam document}
From scattered references to M\=adhava found in diverse manuscripts, historians of Kerala mathematics have pieced together some bits of information. There is a document written in Malayalam script attached to the end of a Malayalam commentary of {\em S\=uryasiddh\=anta} preserved in the Oriental Institute, Baroda (document No. 9886). The document is undated and is composed in a mixture of Malayalam and Sanskrit languages. It gives, in the form of a list, information about a continuous line of scholars from Govinda Bha\d t\d tatiri (born 1237 CE) to Acyuta Pi\d s\=ara\d ti (died 1621 CE). The ninth item in the list gives information about M\=adhava. The item is reproduced below in full (see \cite{SarmaDirect}): 
\begin{quote}
``{\em
M\=adhavan Ila\~n\~nippa\d l\d li Empr\=an Ve\d nv\=aroham, Muh\=urttam\=adhav\=\i yam,\\ Pra\'snam\=adhav\=\i yam mutal\=aya grantha\.n\.na\d  lkkum, Aga\d nita-pa\~nc\=a\.ngattinum\\ kartt\=a\-v\=akunnu.}''\footnote{The words  in this quote appear in a different order in a quote of the same in K. V. Sarma's {\em A History of Kerala School of Hindu Astronomy} (see \cite{Sarma1972} p. 51). In this work, the quote appears thus: ``{\em M\=adhavan v\=e\d nv\=ar\=oh\=ad\=\i n\=a\d m kart\=a ... M\=adhavan Ilaññippa\d l\d li Empr\=an}.''}

This may be translated as follows:
``M\=adhavan Ila\~n\~nippa\d l\d li Empr\=an, author of the books {\em Ve\d nv\=aroham,  Muh\=urttam\=adhav\=\i yam, Pra\'snam\=adhav\=\i yam} etc. and of {\em Aga\d nita-pa\~nc\=a\.nga}.'' 
\end{quote}

It is known for sure that our hero Sa\.ngamagr\=ama M\=adhava is the author of {\em Ve\d nv\=a\-roham}.  It is also known for sure that the author of other two books, {\em Muh\=urttam\=adhav\=\i yam} and {\em Pra\'snam\=adhav\=\i yam}, mentioned in the quoted text  is a different person known as 
Vidy\=am\=adhava\footnote{Vidy\=am\=adhava's {\em Muh\=urttam\=adhav\=\i yam} has been hugely popular in Kerala. It has attracted at least six commentaries by Keralite astronomers two of which are in Sanskrit and the rest in Malayalam. Even though Vidy\=am\=adhava lived in a region far away form the geographical boundaries of present-day Kerala, K. V. Sarma considered him as an illustrious member of the Kerala school of astronomy and mathematics. See \cite{Sarma1972} p. 48.}
 hailing from Gu\d navati, a village near Gokar\d na on the Arabian coast and who is known to have flourished around 1350 CE (see \cite{Sarma1972} p. 48). As per the document, M\=adhava is the author of a fourth work titled {\em Aga\d nita-pa\~nc\=a\.nga} and there are some indirect indications that suggest that the author of this work is Sa\.ngamagr\=ama M\=adhava (see p. \pageref{Aganita}). Vidyam\=adhava, author of {\em Muh\=urttam\=adhav\=\i yam} and {\em Pra\'sna\-m\=adhav\=\i yam}, is also known to be a Tulu Brahmin ({\em Empr\=an}). So the epithet {\em Empr\=an} fits into the name of Vidyam\=adhava and whether it fits into the name of our hero is a matter still under  scholarly scrutiny (see \cite{Sundareswaran2019}). 

Let us examine the meanings of the  terms   ``{\em Empr\=an}'' and   ``{\em Ilaññippa\d l\d li}''  in a little more detail. 
\subsection{``{\em Empr\=an}''}
Let us examine the term ``{\em Empr\=an}'' first. It should be noted that this term actually appears as an epithet to the name ``M\=adhavan''. This epithet ``{\em Empr\=an}'' in the name is a reference to a certain community he would have belonged to. The word ``{\em Empr\=an}''\index{Empr\=an} is a shortened form of the word ``{\em Empr\=antiri}"\index{Empr\=antiri} and the word {\em Empr\=antiri} refers to a member of the community of Br\=ahmins, considered somewhat inferior in status to the Namp\=utiris (see \cite{Plofker2009} p. 218), who have migrated from Tulu Nadu\index{Tulu} to Kerala and settled in Kerala.\footnote{``Empr\=antiri or Empr\=an is a Malayalam name for Tulu Br\=ahmans settled in Malabar. They speak both Tulu and Malayalam. Some of them call themselves Namb\=udiris, but they never intermarry with that  class. Empr\=antiri is a name for Tulu Shivalli Br\=ahmans.'' Edgar Thurston (assisted by K. Rangachari), {\em Castes and Tribes of Southern India}, Volume II, Government Press Madras, 1909, p. 209.} Tulu Nadu is not an administrative unit in modern India and is   generally identified as the region in the southwestern coast of India consisting of the Dakshina Kannada and Udupi districts of Karnataka state and Kasaragod district of Kerala state. Most members of the community speak both Tulu and  Malayalam languages.

%
%
\subsection{``{\em Ilaññippa\d l\d li}''}
\index{Ilaññippa\d l\d li}
The term
``{\em Ilaññippa\d l\d li}"\footnote{Ulloor S. Paramesvara Aiyer, in his history of Kerala literature, has referred to M\=adhava's house name as {\em Iri\~n\~n\=a\d tappa\d l\d li}.}
 has been identified as a reference to the house name of M\=adhava. This is corroborated by M\=adhava himself. In his short work on the moon's positions titled {\em Ve\d nv\=aroha}, M\=adhava says that he was born in a house named {\em baku\d l\=adhi\d s\d thita $\ldots$ vih\=ara} (see \cite{Sphuta} p. 12): 
\begin{quote}
{\dn 
b\7{k}\30FwAEDE\3A4wt(v\?n EvhAro yo Evf\?\309wyt\?.\\
\9{g}hnAmEn so\35Fwy\2 -yAE\3E0wjnAmEn mADv,..
}\\[1mm]
{\em 
baku\d l\=adhi\d s\d thitatvena vih\=aro yo vi\'se\d syate $\scriptstyle\vert$\\
g\d rhan\=am\=ani soya\d m sy\=annijan\=amani  m\=adhava\d h $\scriptstyle\vert\vert$
}
\end{quote}

The Sanskrit term ``{\em baku\d l\=adhi\d s\d thita $\ldots$ vih\=ara}'' is clearly a Sanskritisation of {\em Ilaññi\-ppa\d l\d li}. The word ``{\em ilaññi}'' (or, ``ira\~n\~ni'') is the Malayalam name of the evergreen tree {\em Mimusops elengi} and the Sanskrit name for the same is {\em baku\d la}. The word ``{\em pa\d l\d li}'' had been in use to refer to a Buddhist retreat house and the word continued to be used as suffixes to names of houses and places in Kerala even after the disappearance of Budhism from Kerala. The Sanskrit equivalent of {\em pa\d l\d li} is {\em vih\=ara}. Thus {\em baku\d l\=adhi\d s\d thita} can be translated as ``occupied or inhabited by {\em baku\d la}". Acyuta Pi\d s\=ara\d ti (c.
 1550--1621) in his Malayalam commentary on {\em Ve\d nv\=aroha} has given the same interpretation. According to Pi\d s\=ara\d ti, {\em baku\d la\d m} is {\em ira\~n\~ni} and {\em vih\=ara\d m} is {\em pa\d l\d li} and so the house name is {\em ira\~n\~ni ninna pa\d l\d li}, that is, ``occupied or inhabited by {\em ira\~n\~ni}''. Incidentally, the words {\em ilaññi} and {\em ira\~n\~ni} are synonymous in Malayalam. 
K. V. Sarma has tried to identify the house name {\em ira\~n\~ni ninna pa\d l\d li} with a currently existing Namp\=utiri house with name {\em Iri\~n\~n\=alappa\d l\d li}. He further surmises that the house name {\em ira\~n\~ni ninna pa\d l\d li}, in course of time, might have got transformed into  {\em Iriññanava\d l\d li} or {\em Iriññārappa\d l\d li} and then into {\em Iri\~n\~n\=alappa\d l\d li} (see  \cite{SarmaMathrubhumi}, \cite{Sphuta}). But, according to P. P. Divakaran (see \cite{Divakaran2008} p. 286), this identification is far fetched because these names have neither phonetic similarity nor semantic equivalence to the word ``{\em Ilaññippa\d l\d li}".

\subsection{``{\em Sa\.ngamagr\=ama}''}

N\=\i laka\d n\d tha Somay\=aj\=\i's  {\em \=Aryabha\d t\=\i ya-bh\=a\d sya} is one of the few texts which extensively  discusses M\=adhava's contributions to mathematics. 
In this {\em Bh\=a\d sya}, Somay\=aj\=\i\ 
has at several places  
 referred to M\=adhava as ``M\=adhava born in Sa\.nga\-magr\=ama''. For example,  while discussing circumferences of circles, Somay\=aj\=\i\  invoked M\=adhava's name thus (see \cite{Aryabhatiya1930} p. 42, \cite{Divakaran2008} p. 270):
\begin{quote}
{\dn 
s\3BDwmg\5Amjo mADv, \7{p}nr(yAs\3E0wA\2 pErEDs\2HyA\7{m}\3C4wvA\qq{n}{\rs --\re} 
}
\\[1mm]
{\em sa\.ngamagr\=amajo m\=adhava\d h punaraty\=asann\=a\d m paridhisa\d mkhy\=amuktav\=an}
\\[1mm] ``M\=adhava born in Sa\d ngamagr\=ama ({\em Sa\.ngamar\=ajagjo}) has given an extremely
close ({\em aty\=asanna}) value of the circumference --''
\end{quote}
Also, in another context, Somay\=aj\=\i\ has referred to M\=adhava thus (see \cite{Aryabhatiya1930} p. 58): 
\begin{quote}
{\dn 
\7{p}n-tE\392wqy\2 vs\306wtEtlk\2 s\2gmg\5AmjmADvEnEm\0t\2 p\38Dw\2 c \7{\399w}t\qq{m}.
}
\\[1mm]
{\em punastadvi\d saya\d m vasantatilaka\d m sa\d mgamagr\=amajam\=adhavanirmita\d m padya\d m \\ ca \'sruta\d m $\scriptstyle |$}\\[1mm]
``Again, in this matter, a verse in {\em vasanthathilaka\d m} (meter) composed by M\=adhava born in Sa\d ngamagr\=ama is well known.''
\end{quote}

The commonly accepted  view among most historians and scholars, including K. V. Sarma, is that Sa\.ngama\-gr\=ama\index{Sa\.ngamagr\=ama} is the town of Irinjalakkuda some 70 km south of the Nila river\index{Nila river} and about 70 km south of Cochin (see, for example, \cite{Sarma1958}).\index{Cochin} 
It seems that there is not much concrete ground for this belief except perhaps the fact that the presiding deity of an early medieval temple in the town, the Koodalmanikyam Temple, is worshipped as Sangameswara\index{Sangameswara} meaning the Lord of the Sa\.ngama and so Sa\.ngamagr\=ama can be interpreted as the village of Sa\.ngamesvara (see \cite{Ulloor} p. 98 and \cite{SarmaMathrubhumi}). 

But there are several places in Karnataka with {\em sa\.ngama} or its Dravidian equivalent {\em kū\d dala} in their names and with a temple dedicated to Sa\.ngamesvara, the lord of the confluence. (Kudalasangama in Bagalkote district in Karnataka State, which lies at the point of confluence of the rivers Krishna and Malaprabha, is one such place with a celebrated temple dedicated to the Lord of the Sa\.ngama.) This should also be read in conjunction with the fact that M\=adhava had been referred to as an {\em Empr\=an} and, as we have already noted earlier, an {\em Empr\=an} is a member of a community whose members had migrated from Tulu Nadu to Kerala and settled in Kerala.

\begin{figure}[t]
\centering
\includegraphics[height=4.7cm]{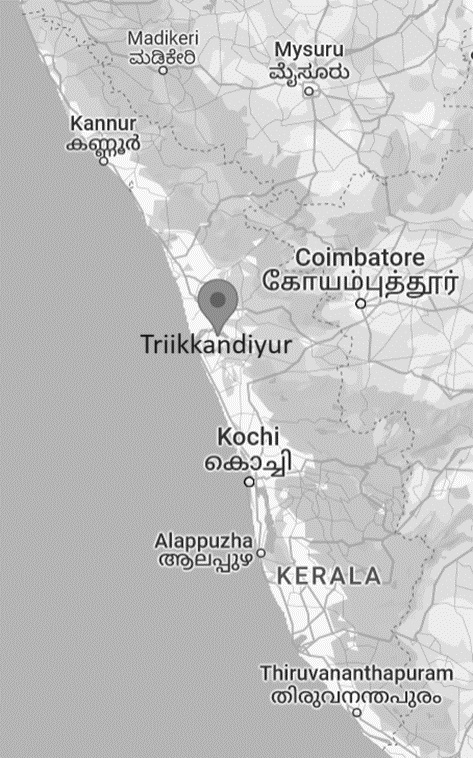}\quad 
\includegraphics[height=4.7cm]{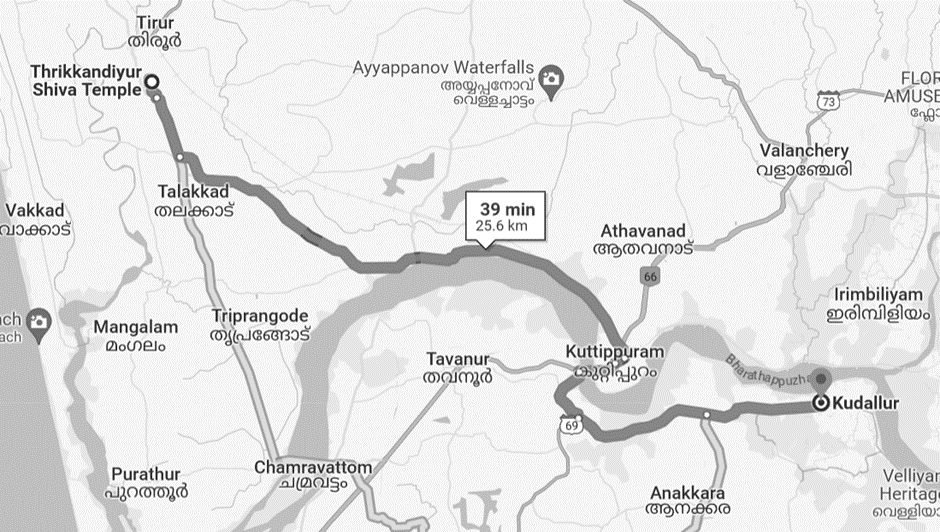}
\caption{Google Map showing the road connecting Thrikkandiyur and Kudallur. Many prominent members of the Kerala School of Astronomy and Mathematics, including N\=\i laka\d n\d tha Somay\=aj\=\i, flourished near Thrikkandiyur.}\label{KudallurMap}
\end{figure}

Interestingly, there is a small town on the southern banks of the Nila river, around 10 km upstream from Tirunnavaya and about 25 km by road from Thrikka\-ndiyur, called K\=u\d tall\=ur (sometimes transliterated as Kudallur) (see Fig.\ref{KudallurMap}).\index{Kudallur@K\=u\d tall\=ur} Thrikkandiyur has a prominent place in the history of Kerala school of astronomy and mathematics. Many prominent members of the school, including N\=\i laka\d n\d tha Somay\=aj\=\i\  and Acyuta Pi\d s\=ara\d ti, flourished near Thrikkandiyur. The exact literal Sanskrit translation of this place name is Sa\.ngamagr\=ama: ``({\em k\=u\d tal})'' in Malayalam means a confluence (which in Sanskrit is {\em sa\.ngama}) and ``{\em \=ur}'' means a village (which in Sanskrit is {\em gr\=ama}). Also the place is at the confluence of the Nila river and its most important tributary, namely, the Kunti river. (There is no confluence of rivers near Irinjalakuda.) 

Incidentally there is a still existing Namp\=utiri family by name K\=u\d tall\=ur house (locally called K\=u\d tall\=ur Mana) a few kms away from the K\=u\d tall\=ur village. (In fact this is not a new observation. As early as 1945, Dr Kunjunni Raja had made a similar observation (see \cite{SarmaThesis} p. 1268).) 
The family has its origins in K\=u\d tall\=ur village itself. In this context, it is interesting and instructive to observe that, there is a Sanskrit drama titled {\em Kamalin\=\i-ka\d laha\d msa} composed around 17th-18th centuries CE which contains a reference to ``{\em sa\.ngamagr\=ama g\d rha}'', that is a ``{\em sa\.ngamagr\=ama} house'', that is ``{\em K\=u\d tall\=ur} house'':
\begin{quote}
{\dn 
aE-t k\?rl\?\7{q} s\3BDwmg\5Am\2 nAm \9{g}h\qq{m}.
}%

{\em
asti kerale\d su sa\.ngamagr\=ama\d m n\=ama g\d rha\d m $\scriptstyle\vert$
}
\end{quote}
It is known that the author of the drama is one N\=\i laka\d n\d tha belonging to K\=u\d tall\=ur house. In the drama, the reference to {\em sa\.ngamagr\=ama g\d rha} appears in a prologue as part of the self-description of the author (see \cite{Iyer1971} p. 605,  \cite{Sundareswaran2019}). 
This corroborates the fact that the term {\em sa\d ngamagr\=ama} had been in use as a Sanskritisation of the  name of the K\=u\d tall\=ur house, and probably by extension, of the K\=u\d tall\=ur village also.

Sanskritisation of proper nouns, especially place names and house names, was a common practice among Kerala authors of Sanskrit works. We have already seen that M\=adhava himself has Sanskritised his house name {\em Ila\~n\~nippa\d l\d li} as {\em Bakul\=adhi\d s\d thita-vih\=ara}. 
N\=ila\-ka\d n\d tha Somay\=aj\=\i\ (1444--1545), a prominent figure of the Kerala School of astronomy and mathematics, refers to his house name as {\em Kerala-sadgr\=ama} which is a Sanskritisation of his Malayalam house name {\em K\=e\d lall\=ur} (taking it as a corrupted form of {\em Kerala-nal-\=ur}),  and the village where he lived as \'Sr\=\i\-kund\=apura or \'Sr\=\i \-kund\=a\-gr\=ama which is a Sanskritisation of T\d rkka\d n\d tiy\=ur, the local name of the village. The name ``Jye\d s\d thadeva'' of the author of the celebrated {\em Yuktibh\=a\d s\=a} is, most probably, a Sanskritisation of the his real personal name (see \cite{Sarma1972} for more such instances of Sanskritisation).

For many generations the K\=u\d tall\=ur family hosted a great Gurukulam. Perhaps it is not a coincidence that  the only available manuscript of {\em Sphu\d tacandr\=apti}, a book authored by M\=adhava, was obtained from the manuscript collection of K\=u\d tall\=ur Mana (see \cite{Sphuta} p. 8). M\=adhava might have had some association with K\=u\d all\=ur house. K. V. Sarma has explicitly used a similar argument to ascertain the identity of another work. Based on the fact that the manuscript of a commentary on a certain work titled {\em K\=alad\=\i pa} has been found only in the palace of {\em Punnatt\=urk\=o\d t\d ta}, Sarma has concluded that the author of the commentary is likely to be a member of the {\em Punnatt\=urk\=o\d t\d ta}\index{Punnatt\=urk\=o\d t\d ta} household (see \cite{Sarma1972} p.70). 

\section{M\=adhava's times}

Regarding M\=adhava's times there are only a few indirect evidences to determine the period during which M\=adhava lived. In spite of very many diligent efforts by historians of Kerala mathematics, nothing concrete has emerged to pinpoint the exact period during which M\=adhava flourished.  In fact, in the astronomical and mathematical literature of Kerala, there is a deafening silence about M\=adhava the person, about his family, about his parents, about his siblings and even about his teachers. But, fortunately in this literature there are several references to one of his highly respected student about whom some concrete information is available, and from this information we can deduce approximately the period in which M\=adhava lived and worked. Also, M\=adhava has inadvertently left some clues in his works which indicate  the historical period in which he flourished. 

\subsection{We know his student!}

Even though, at present, we do not know anything about the the parents or siblings of M\=adhava, we do know much about one of M\=adhava's illustrious pupil, namely, Va\d{t}a\'sre\d{n}i Parame\'svara. N\=\i laka\d n\d tha Somay\=aj\=\i\ in his {\em \=Aryabha\d t\=\i ya-bh\=a\d sya}, while talking about his own teachers, refers to Parame\'svara and his son D\=amodara as  his great teachers and and also mentions that Parame\'svara had mastered astronomy and mathematics from, among others, M\=adhava the {\em Golavid}\index{Golavid@{\em Golavid}} (see \cite{SarmaDirect}). There are some clear references associating certain specific years with the life of Parame\'svara. Parame\-\'svara has stated  in his {\em d\d rgga\d nia} that he completed the work in \'Saka 1353, equivalently 1431 CE and his {\em Golad\=\i  pik\=a} in \'Saka 1365 (1443 CE). Parame\'svara is reputed to have composed his {\em d\d rgga\d nia} system based on his observations of the skies extending over a period of fifty-five years. This would mean that he had commenced his investigations by about 1375 CE. This coupled with the fact that N\=\i laka\d n\d tha Somay\=aj\=\i\  (born 1443 CE) was a pupil of Parame\'svara would make the probable period of Parame\'svara around 1360--1460. 

\subsection{Hidden references in works}

Inadvertently though, M\=adhava has left some clues in his works regarding the times in which he lived and worked. These clues properly interpreted
 shed some light on M\=adhava's times.
M\=adhava has been identified as the author of a short treatise called {\em Sphu\d ta\-candr\=apti} and also of another work titled {\em Ve\d nv\=aroha} which is an elaboration of {\em Sphu\d t\-acandr\=apti}. These works basically describe an ingenious method to compute the true positions of the Moon at intervals of about half an hour for various days in an anomalistic cycle. 

In these works, for the calculation of the position of the Mean Sun\footnote{The {\em Mean Sun}\index{Mean Sun} is an imaginary sun moving along the celestial equator at a constant rate and completing its annual course in the same time as the real sun takes to move round the ecliptic at varying rates.}, we are asked to subtract from the current {\em kali} day\index{kali day@{{\em kali} day}}, a duration of 1,502,008 days and also the duration of 5180 anomalistic cycles of the Moon. If there be further days, the number of such days has to be multiplied by the mean daily motion and added. It has been estimated that the day corresponding to the 1,502,008-th  {\em kali}-day followed by 5180 anomalistic cycles of the Moon is 10 March 1402 CE (see 
  \cite{Hari2003}, \cite{Sarma1958}, \cite{Sphuta}). Since there is no astronomical significance to this date, it is reasonable to assume that this date would not be far away from the date of completion of {\em Sphu\d tacandr\=apti} or {\em Ve\d nv\=aroha}. This suggests that M\=adhava should have lived around 1400 CE.

In {\em Aga\d nita}, another work identified as M\=adhava's, the author has given the {\em \'sodhy\=a\-bda}-s (`deductive years') for computing the positions of various planets. The `deductive years' represent the largest number of years that could be cut off for the different planets at the time when the work was written. The largest deductive year mentioned in the work corresponds to 1418 CE. (M\=adhava has specified the deductive years as years in the {\em \'S\=aka} era.) From this, one can reasonably conclude that the date of completion of {\em Aga\d nita} is just ahead of this the year 1418 CE (see \cite{Sphuta}).

\subsection{So M\=adhava lived during . . .}

Based on these evidences it appears reasonable to assume that M\=adhva flourished around the period 1340--1425 CE (see \cite{Sarma1958}).
\subsection{A little more about M\=adhava's student}
\index{Va\d{t}a\'sre\d{n}i Parame\'svara}
\index{Parame\'svara, Va\d{t}a\'sre\d{n}i}
\index{drigganita@{\em d\d{r}gga\d{n}ita}}
\index{Haridatta}
\index{Alathiyur}
Va\d{t}a\'sre\d{n}i Parame\'svara (Va\d{ta}\'s\'seri Parame\'svaran Namp\=utiri) (c. 1360--1460), who introduced the {\em d\d rgga\d{n}ita} system for astronomical computations and who is considered as one of greatest astronomers of the Kerala school, was an illustrious pupil of Sa\.ngama\-gr\=ama M\=adhava. The {\em d\d rigga\d{n}ita} system propounded by Parame\'svara was a revision of the {\em parahita} \index{parahita@{\em parahita}} system introduced by Haridatta \index{Haridatta} in the year 683 CE. No new methodology was introduced as part of the {\em d\d rigga\d{n}ita} system. Instead, new multipliers and divisors were given for the computation of the {\em kali} days and for the calculation of the mean positions of the planets. Revised values are given for the positions of planets at zero {\em kali}. Also the values of the sines of arc of anomaly ({\em manda-jy\=a}) and of commutation ({\em s\=\i ghra-jy\=a}) are revised and are given for intervals of 6 degrees. Parame\'svara was also a prolific writer on matters relating to astronomy. At least 25 manuscripts have been identified as being authored by Parame\'svara.

 Parame\'svara's family name was Vata\'s\'seri\index{Vata\'s\'seri} also written as Va\d{t}a\'sre\d{n}i and his family resided in the village of Alathiyur (Sanskritised as A\'svatthagr\=ama) in Tirur, Kerala. Ala\-thi\-yur is situated on the northern bank of the  Nila river at its mouth.(See \cite{SarmaThesis} Chapter XIII {\em D\d rgga\d nita of Parame\'svara}, pp. 1049--1094.) 
\section{M\=adhava's works}
The following is a list of works whose authorship has been attributed to M\=adhava.  As has  been noted by several historians of Kerala mathematics, this list is unlikely to be an exhaustive list of all the works of M\=adhava because only a fraction of all the astronomical and mathematical manuscripts from Kerala has been studied, edited and translated and it is quite possible that many other  great works of M\=adhava may be lying buried deep in the mountains of unexplored manuscripts.   Also none of the works listed below (except perhaps {\em Mah\=ajy\=anayanaprak\=ara} which is of doubtful authorship) contains details of the most celebrated mathematical achievements of M\=adhava.
\begin{enumerate}
\item
{\em Ve\d{n}v\=aroha} (Ascending the bamboo):\index{Ve\d{n}v\=aroha@{\em Ve\d{n}v\=aroha}} This is a work in 74 verses describing methods for the computation of the true positions of the Moon at intervals of 36 minutes for various days in an anomalistic cycle. This work is an elaboration of an earlier and shorter work of Ma\=dhava himself titled {\em Sphu\d tacandr\=apti}. {\em Ve\d{n}v\=aroha} is the most popular astronomical work of M\=adhava. The novelty and ingenuity of the method caused several works of the same genre to be composed by later astronomers. Thus we have {\em Ve\d nv\=aro\-h\=a\d s\d taka} by Putumana Somay\=aj\=\i\ and and {\em D\d rg-ve\d nv\=arohakriy\=a} of epoch 1695 of anonymous authorship. 

K. V. Sarma has brought out a critical edition of {\em Ve\d nv\=aroha} with Acyuta Pi\d s\=ara\d ti's commentary in Malayalam, and the same is available online for free download (see \cite{SarmaThesis} Chapter XVI {\em Ve\d nv\=aroha}, pp. 1259--1316). Acyuta Pi\d s\=ara\d ti (c. 1550--1621) was a student of Jye\d s\d thadeva, the author of the celebrated {\em Yuktibh\=a\d s\=a}, and authored several treatises on astronomy.
\item
{\em Sphu\d tacandr\=apti} (Computation of true Moon):\index{Sphu\d tacandr\=apti@{\em Sphu\d tacandr\=apti}}  
 This is a shorter version of M\=a\-dhava's {\em Ve\d nv\=aroha} and discusses nearly all topics that have been included in {\em Ve\d nv\=a\-roha}. 
 
 A critical edition of {\em Sphu\d tacandr\=apti} along with an English translation has been brought out by K. V. Sarma and the same is available for free download from the internet (see \cite{Sphuta}).
\item
{\em Candrav\=aky\=ani}\index{Candrav\=aky\=ani@{\em Candrav\=aky\=ani}} (Moon-sentences): 
This is sometimes considered as an independent work of M\=adhava. 
This is a collection of 248 numbers, arranged in the form of a list, related to the motion of the Moon in its orbit around the Earth. These numbers are couched in the {\em katapayadi} system and so appear like a list of words, or phrases or short sentences written in Sanskrit. M\=adhava's list is a refinement of an earlier list, also called {\em Candrav\=aky\=ani}, attributed to the legendary Vararuci. Vararuci's values are correct only to the arc-minute whereas M\=adhava's values are correct to the arc-second. 

The text of M\=adhava's {\em Candrav\=aky\=ani} has been included as appendices in the modern publications of {\em Ve\d{n}v\=aroha} and {\em Sphu\d tacandr\=apti} (see \cite{SarmaThesis} Chapter XVI {\em Ve\d nv\=a\-roha}, Appendix I, pp. 1309--1312 and \cite{Sphuta} Appendix I, pp. 46--57). 
\item
{\em Mah\=ajy\=anayana-prak\=ara} (Method for the computation of the great sine):\index{Mah\=ajy\=anayanaprak\=ara@{\em Mah\=ajy\=anayanaprak\=ara}} K. V. Sarma attributed\footnote{David Pingree, an American historian of mathematics in the ancient world, expressed doubts about this attribution (see \cite{Pingree1991}). However, Pingree agreed that the author should be from `M\=adhava school'.} the authorship of this work  to M\=adhava. Probably, this is the only surviving work  by M\=adhava   dealing with the power series expansions of the sine and cosine functions. This work contains novel theorems and computational methods developed by M\=adhava and used by later astronomers and mathematicians of the Kerala school. 
The manuscript of the work has three sections, namely, the explanation of the
series, the method to derive the numerical sine values and the
derivation for both the sine and cosine series.

{\em Mah\=ajy\=anayanaprak\=ara} has been critically edited and translated into English by David Pingree  (see \cite{Pingree1991}). More details about this work are available in \cite{Rajeswari}.
\item
{\em Madhym\=anayanaprak\=ara}\index{Madhym\=anayanaprak\=ara@{\em Madhym\=anayanaprak\=ara}} (The procedure for obtaining the mean):
This short work describes a procedure for the computation of the {\em madhyamagraha}
(mean longitude of a planet). Like the Mean Sun, a ``mean planet'' is a hypothetical planet that coincides with a real planet when the real planet is at perihelion and that moves in an orbit at a constant velocity equal to the mean velocity of the real planet. The longitude of the mean planet is known as the ``mean longitude'' or ``{\em madhyamagraha}''. 

This short tract starts with a description of the geometrical construction involved in the {\em manda-sa\d msk\=ara} and then  proceeds to describe in detail the concept of {\em vipar\=\i ta\-kar\d na}. The treatise then presents a detailed discussion on its application to simplify the 
the computation of {\em avi\'si\d ta-manda-kar\d na}. It concludes with a description of the procedure to compute {\em madhyamagraha}. 

 The work is essentially a very detailed commentary in prose on just two verses: one beginning with {\em vist\d rti} and the other beginning with {\em arkendu}. Both verses appear in N\=\i laka\d n\d tha Somay\=aj\=\i's {\em Tantrasa\.ngraha}, but Somayaj\=\i\ has ascribed the authorship of the verses to M\=adhava. The ascription of the prose commentary to M\=adhava is due to K. V. Sarma (see \cite{Madhayama}).

An edited version of the text along with English translation
 and explanations of  the technical content  using modern mathematical notations has been published as a research paper by U. K. V. Sarma and others
(see \cite{Madhayama}).
\item
{\em Aga\d{n}ita}/{\em Aga\d{n}ita-grahac\=ara}/{\em Aga\d{n}ita-pa\~nc\=a\.nga}\index{Aga\d nita@{\em Aga\d nita}}\index{Aga\d{n}ita-grahac\=ara@{\em Aga\d{n}ita-grahac\=ara}}:\label{Aganita}  According to an undated astronomical palm-leaf document written in Malayalam script preserved in the Oriental Institute, Baroda, M\=adhava is the author of a work titled {\em Aga\d{n}ita-pa\~nc\=a\.nga}. A treatise of anonymous authorship titled {\em Aga\d{n}ita-grahac\=ara} is available in manu\-script form\footnote{For details of the manuscript, see Item No. 628 {\em Aga\d nitagrahac\=ara} (p. 1305) in {\em A Descriptive Catalogue of Sanskrit Manuscripts in the Curator's Office Library, Trivandrum Vol. IV} published in 1939. See \cite{Catalogue}.} and Putumana Somay\=aj\=\i's {\em Kara\d napaddhati} contains verses quoted from this treatise. Probably, both {\em Aga\d{n}ita-grahac\=ara} and {\em Aga\d{n}ita-pa\~nc\=a\.nga} refer to the same work (see \cite{Sarma1972} p. 51, \cite{SarmaDirect}).
{\em Aga\d{n}ita-grahac\=ara} contains the longitudes of the planets for long cycles of years in the form of tables.
 
\item
{\em Lagnaprakara\d{n}a} (Treatise on the ascendant)\index{Lagnaprakara\d{n}a@{\em Lagnaprakara\d{n}a}}: This is an astronomical text devoted to the determination of {\em lagna} (ascendant). The term {\em lagna} commonly refers to {\em udaya lagna} which is  the zodiac sign rising on the east at any given instant. The work consists of 139 verses divided into eight chapters and gives the exact values of {\em lagna} unlike earlier texts (see \cite{Kolachana}). In the first chapter of the {\em Lagnaprakara\d na}, M\=adhava
discusses several procedures  to
determine astronomical quantities such as the {\em pr\=a\d na-kal\=antara} (difference between the longitude and right ascension), {\em cara} (ascensional difference of a body),
and {\em k\=alalagna} (the time interval between the rise of the
vernal equinox and a desired later instant).  From the second chapter onwards, the text describes several techniques for precisely
determining the {\em udayalagna} (see \cite{Kolachana2}).

More information about {\em Lagnaprakara\d na} is available in a doctoral thesis titled ``A critical Study of M\=adhava's {\em Lagnaprakara\d na}'' authored by Aditya Kolachana, prepared under the supervision of K. Ramasubramanian   and submitted to IIT Bombay in 
2018. Since then Kolachana, Ramasubramanian and others have published a series of articles expounding in detail the various mathematical and astronomical aspects of {\em Lagnaprakara\d na} (see, for example, \cite{Lagnaprakarana08}, \cite{Lagnaprakarana03}).
\item
{\em Golav\=ada} (Spherics):\index{Golav\=ada@{\em Golav\=ada}} Later writers have often referred to M\=adhava as {\em Go\d lavid} (master of spherics). Researchers have spotted a reference to the existence of a text titled {\em Golav\=ada} authored by M\=adhava in a palm-leaf manuscript of {\em Ve\d n\=aroha}.  It is believed that this work, copies of which are still not extant, earned M\=adhava the appellation {\em Go\d lavid} (see \cite{Sarma1972} p. 52).
\item
{\em . . .  ?} : According to K. V. Sarma ``it seems quite possible that M\=adhava had composed a comprehensive treatise on astronomy and mathematics, which yet remains to be identified and which may be supposed to contain the numerous single and groups of verses enunciating computational procedures, theorems and formulas which are quoted as M\=adhava's by later writers (see \cite{Sarma1972} p. 52).''
\end{enumerate}
\section{M\=adhava's contributions to mathematics in a nutshell}
As a curtain raiser, we present here some of the most important contributions of M\=adhava to mathematics. All of these results, their proofs and other  results associated with them have  been discussed in much fuller  detail in the later chapters of this book. 
\begin{enumerate}
\item
{\bf M\=adhava-Leibniz series}
\index{M\=adhava-Leibniz series}

M\=adhava stated and proved the following infinite series expression for the mathematical constant $\pi$:
$$
\frac{\pi}{4} = 1 -\frac{1}{3}+\frac{1}{5}-\frac{1}{7}+ \cdots
$$
\item
{\bf Correction terms to M\=adhava-Leibniz series}
\index{correction terms}

M\=adhava gave the following correction terms ($F_1(n)$, $F_2(n)$ and $F_3(n)$) to the M\=adha\-va-Leibniz series which help  us obtain better approximations to $\pi$ much more quickly: 
$$
\frac{\pi}{4} \approx 1 -\frac{1}{3}+\frac{1}{5}-  \cdots +(-1)^{n-1}\frac{1}{2n-1} + (-1)^nF_i(n)
$$
where
\begin{align*}
F_1(n)& = \frac{1}{4n},\\
F_2(n)& = \frac{n}{4n^2+1},\\
F_3(n)& = \frac{n^2+1}{n(4n^2+5)}.
\end{align*} 
\item
{\bf Many more infinite series expressions for $\pi$}

In addition to the M\=adhava-Leibniz series, M\=adhava gave several more infinite series expressions for $\pi$ like the following:
\begin{enumerate}
\item
$\dfrac{\pi}{4}  = \dfrac{3}{4} + \dfrac{1}{3^3-3}-\dfrac{1}{5^3-5}+\dfrac{1}{7^3-7} - \cdots$
\item
$\dfrac{\pi}{8}  = \dfrac{1}{2^2-1} + \dfrac{1}{6^2-1} +\dfrac{1}{10^2-1}+\cdots$
\item
$\dfrac{\pi}{4}  = \dfrac{4}{1^5 + 4\cdot 1} - \dfrac{4}{3^5+4\cdot 3}+\dfrac{4}{5^5+4\cdot 5} - \dfrac{4}{7^5+4\cdot 7}+ \cdots$
\item
$\dfrac{\pi}{4}  = \dfrac{1}{2} +\dfrac{1}{2^2-1} -\dfrac{1}{4^2-1}+\dfrac{1}{6^2-1} - \cdots$
\end{enumerate}
\item
{\bf The M\=adhava-Gregory series}
\index{M\=adhava-Gregory series}
\index{Gregory, James}

M\=adhava stated and proved the following infinite series expression for the arctangent   ($\tan^{-1}$) function  
$$
\arctan x = x -\frac{x^3}{3}+\frac{x^5}{5} -\frac{x^7}{7}+\cdots
$$
and the following infinite series expression for $\pi$ derived from the series for $\arctan x$:
$$
\pi=\sqrt{12}\Big( 1 - \frac{1}{3\cdot 3}+
\frac{1}{5\cdot 3^2} - \frac{1}{7\cdot 3^3} + \cdots \Big).
$$
\item
{\bf Numerical value of $\pi$ correct to the 10-th decimal place}

The discovery of the following approximate value of $\pi$, which is correct to the 10-th decimal place, has been attributed to M\=adhava:
$$
\pi \approx  \frac{2,827,433,388,233}{9\times 10^{11}}.
$$
\item
{\bf The M\=adhava-Newton series}

M\=adhava stated and proved the following power series expressions for $\sin \theta$, $\cos\theta$ and $\sin^2\theta$:
\begin{enumerate}
\item
$\sin \theta 
= \theta-\dfrac{\theta^3}{3!}+\dfrac{\theta^5}{5!}- \dfrac{\theta^7}{7!} +\cdots$
\item
$\cos \theta 
= 1 - \dfrac{\theta^2}{2!} + \dfrac{\theta^4}{4!} - \dfrac{\theta^6}{6!}+ \cdots $
\item
$\sin^2\theta
 = 
\theta^2 - \dfrac{\theta^4}{(2^2 -\tfrac{2}{2})}
+
\dfrac{\theta^6}{(2^2 -\tfrac{2}{2})(3^2 -\tfrac{3}{2})}
- \cdots$
\end{enumerate}
\item
{\bf Algorithms to compute $\sin\theta$ and $\cos\theta$}

M\=adhava developed algorithms to compute the numerical values of $\sin\theta$ and $\cos\theta$ using precomputed polynomial coefficients and a very efficient and completely no\-vel method to compute the values of polynomials. This method of computing values of polynomials was discovered in Europe only in the nineteenth century (see \cite{Krish2010}). 
\item
{\bf Computation of a sine table}
\index{sine table}

M\=adhava constructed a table giving the values of $\sin\theta$ for 
$$
\theta = 3.75^\circ, 7.50^\circ, 11.25^\circ, \ldots, 86.25^\circ, 90^\circ.
$$
 The values are correct to the 7th - 8th decimal places.
\item
{\bf M\=adhava's Taylor series approximations to $\sin\theta$ and $\cos\theta$}
\index{Taylor series}

M\=adhava obtained the following Taylor series approximations to the sine and cosine functions:
\begin{align*}
\sin(u+h) 
& = \sin (u) + h\cos (u) - \frac{h^2}{2}\sin (u),\\
\cos(u+h) 
& = \cos (u) - h\sin (u) +\frac{h^2}{2}\cos (u).
\end{align*}
\item
{\bf M\=adhava's {\em j\=\i ve-paraspara-ny\=aya}}
\index{jive-paraspara-nyaya@{\em j\=\i ve-paraspara-ny\=aya}}

The terminology {\em j\=\i ve-paraspara-ny\=aya}  (which may be translated as ``rule of mutual sines'') refers to the following set of rules which specify the addition and subtraction formulas for sine and cosine functions. 
\begin{enumerate}
\item
$\sin(x+y)  = \sin x \cos y + \cos x \sin y$
\item
$\sin(x-y)  = \sin x \cos y - \cos x \sin y$
\item
$\cos(x+y)  = \cos x \cos y - \sin x \sin y$
\item
$\cos(x-y)  = \cos x \cos y + \sin x \sin y$
\end{enumerate}
Even though these rules had been  stated by earlier Indian mathematicians like Bh\=a\-skara II, in the astronomical and mathematical literature of the Kerala school the discovery of these rules has been attributed to M\=adhava.
\item
{\bf Properties of cyclic quadrilaterals}
\index{cyclic quadrilateral}

Va\d ta\'sre\d ni Parame\'svara, one of M\=adhava's distinguished students, in one of his writings, has stated  a formula to compute the circumradius $R$ of a cyclic quadrilateral in terms its sides $a, b, c, d$:
$$
R = \left(\frac{(ab+cd)(ac+bd)(ad+bc)}{(b+c+d-a)(a+d+c-b)(a+b+d-c)(a+b+c-d)}\right)^{\frac{1}{2}}.
$$
If we consider the fact that M\=adhava was a teacher of Parame\'svara, it could be the case that M\=adhava was the originator of the formula even though there is no substantive evidence to prove this surmise. 
\item
{\bf The {\em Ve\d nv\=aroha} algorithm}

{\em Ve\d nv\=aroha} and {\em Sphu\d tacandr\=apti}, with almost identical content, are the only works of M\=adhava that have survived in full to the current times. Though,  the content of {\em Ve\d nv\=aroha} cannot be construed as dealing with mathematics, it is interesting in itself because it gives a detailed procedure, exactly in the form of an algorithm, to compute the true positions of the Moon at various times of a day.
\end{enumerate}
\section{M\=adhava might have made a mistake}
M\=adhava was great, but he was human too and so knowing that he might have made a silly mistake in a numerical computation only would add more feathers to his greatness. M\=adhava has been quoted by several later authors as saying that the circumference of a circle of diameter $9\times 10^{11}$ is $2,827,433,388,233$. The correct value of the circumference is $2,827,433,388,231$. What is strange about this is not the difference between the two values, but the fact that all the formulas for computing $\pi$ attributed to M\=adhava and his followers yield  values for the circumference of a circle of diameter $9\times 10^{11}$  closer to the true value than the value attributed to M\=adhava. Probably M\=adhava made a mistake in his computations! 
 (See \cite{Krish2024X}.)


\begin{thebibliography}{10}

\bibitem{Divakaran2008}
P.~P. Divakaran.
\newblock {\em The Mathematics of {I}ndia: {C}oncepts, methods and
  Connections}.
\newblock Springer/Hindustan Book Agency, New Delhi, 2008.
\newblock ISBN 978-93-86279-65-1.

\bibitem{Pingree1991}
David Gold and David Pingree.
\newblock A hitherto unknown {S}anskrit work concerning {M}adhava's derivation
  of the power series for sine and cosine.
\newblock {\em Historia Scientiarum}, 42:49--65, 1991.
\newblock Available at \url{https://dl.ndl.go.jp/pid/11023753/1/28}.

\bibitem{Hari2003}
K.~Chandra Hari.
\newblock Computation of true moon by {M}\=adhava of {S}a\.ngamagr\=ama.
\newblock {\em Indian Journal of History of Science}, 38(3):231--253, 2003.
\newblock Available at
  \url{https://web.archive.org/web/20120316083104/http://www.new.dli.ernet.in/rawdataupload/upload/insa/INSA_1/2000c4df_231.pdf}.

\bibitem{Ulloor}
Ulloor S.~Parameswara Iyer.
\newblock {\em {Kerala Sahithya Charithram (History of Kerala Literature) Vol.
  II}}.
\newblock University of Travancore, Thiruvananthapuram, 1954.
\newblock Available at \url{https://books.sayahna.org/ml/pdf/ulloor-vol-2.pdf}.

\bibitem{Kolachana2}
Aditya Kolachana, K.~Mahesh, and K.~Ramasubramanian.
\newblock Precise determination of the ascendant in the {L}agnaprakara\d na -
  {I}.
\newblock {\em Indian Journal of History of Science}, 54(3):304--316, September
  2019.
\newblock Available at
  \url{https://cahc.jainuniversity.ac.in/assets/ijhs/Vol54_3_2019__Art04.pdf}.

\bibitem{Lagnaprakarana08}
Aditya Kolachana, K.~Mahesh, and K.~Ramasubramanian.
\newblock {Precise Determination of the Ascendant in the {\em Lagnaprakara\d na
  - IV}}.
\newblock {\em Indian Journal of History of Science}, 56(2):1--13, July 2021.
\newblock Available at
  \url{https://www.researchgate.net/publication/354592918_Precise_determination_of_the_ascendant_in_the_Lagnaprakarana-IV}.

\bibitem{Kolachana}
Sri Ram~Aditya Kolachana.
\newblock {\em Critical study of {M}\=adhava's Lagnaprakara\d na}.
\newblock 2018.
\newblock (Ph.D. thesis submitted to Indian Institute of Technology Bombay,
  Department of Humanities and Social Science.).

\bibitem{Krish2024X}
V.~N. Krishnachandran.
\newblock {\em On \'Sankara Varman's (correct) and M\=adhava's (incorrect)
  values for the circumferences of circles}.
\newblock Online, May 2024.
\newblock (Preprint) Available at \url{https://arxiv.org/abs/2405.11144}.

\bibitem{Krish2010}
V.~N. Krishnachandran, Reji~C. Joy, and Siji~K. B.
\newblock On sangamagr\=ama m\=adhava's (c.1350 -- c.1425 ce) algorithms for
  the computation of sine and cosine functions.
\newblock In {\em Proceedings of the International Conference on Computational
  Engineering Practices and Technology}. MES College of Engineering,
  Kuttippuram, Kerala, November 2010.

\bibitem{Plofker2009}
Kim Plofker.
\newblock {\em Mathematics in {I}ndia}.
\newblock Princeton University Press, Princeton, 2009.

\bibitem{Rajeswari}
G.~Raja Rajeswari and M.~S. Sriram.
\newblock {Mah\=ajy\=anayanaprak\=ara\d h: Infinite Series for the Sine and
  Cosine Functions in the Kerala Works}.
\newblock In Sita~Sundar Ram and Ramakalyani V., editors, {\em {History and
  Development of Mathematics in India}}, pages 294--306. National Manuscript
  Mission, New Delhi, 2022.
\newblock Available at
  \url{https://namami.gov.in/sites/default/files/book_pdf/History%20and%20Development%20of%20Mathematics%20in%20India.pdf}.

\bibitem{Iyer1971}
Subramania~Iyer S.
\newblock {\em {A critical study of Sanskrit dramas by Kerala authors}}.
\newblock University of Kerala, 1971.
\newblock PhD thesis available at
  \url{https://shodhganga.inflibnet.ac.in/jspui/bitstream/10603/142047/18/18_chapter%2012.pdf}.

\bibitem{SarmaMathrubhumi}
K.~V. Sarma.
\newblock {Sa\.ngamagr\=ama M\=adhavan}.
\newblock {\em Mathrubhumi Weekly}, Nov. 1956.
\newblock (The article is in Malayalam.).

\bibitem{Sarma1958}
K.~V. Sarma.
\newblock Date of {M}\=adhava - {A} little known {I}ndian astronomer.
\newblock {\em Quarterly Journal of the Mythic Society, Bangalore},
  XLIX(3):183--186, 1958.
\newblock Available in {\em Facets of Indian Astronomy} (A Collection of
  Articles of Prof. K. V. Sarma), edited by Prof Siniruddha Dash, Rashtriya
  Sanskrit University, Tirupati, 2009, pp.270–273 which is available at
  \url{https://archive.org/details/facetsofindianastronomyarticlesofsarmak.v.ed.siniruddhadashvenkateswarauniversity_840_N/page/270}.

\bibitem{Sarma1972}
K.~V. Sarma.
\newblock {\em A History of the {K}erala School of {H}indu Astronomy (in
  perspective)}.
\newblock Visheshvaranand Vishvabandhu Institute of Sanskrit and Indological
  Studies, Panjab University, Hoshiarpur, 1972.
\newblock Available at
  \url{https://archive.org/details/KeralaSchoolOfAstronomy}.

\bibitem{Sphuta}
K.~V. Sarma.
\newblock {\em {Computation of the True Moon by M\=adhava of Sangamagr\=ama
  (Critically edited with introduction, translation and notes)}}.
\newblock Visheshvaranand Vishvabandhu Institute of Sanskrit and Indological
  Studies, Panjab University, Hoshiarpur, 1973.
\newblock Available at \url{https://archive.org/details/SphutaChandrapti}.

\bibitem{SarmaThesis}
K.~V. Sarma.
\newblock {\em Contributions to the study of {K}erala school of {H}indu
  astronomy and mathematics}, volume~II.
\newblock 1977.
\newblock Thesis submitted to the Panjab University for the Degree of Doctor of
  Philosophy. Available at \url{
  https://archive.org/details/contributionstothestudyofthekeralaschoolhinduastrnomymathematicssarmak.v.vol2thesis_997_}.

\bibitem{SarmaDirect}
K.~V. Sarma.
\newblock Direct lines of astronomical tradition in {K}erala.
\newblock In Sinirudhha Dash, editor, {\em {Facets of Indian Astronomy}}, pages
  402--405. Rashtriya Sanskrit University, Tirupai, Andhra Pradesh, 2009.
\newblock Available at
  \url{https://archive.org/details/facetsofindianastronomyarticlesofsarmak.v.ed.siniruddhadashvenkateswarauniversity_840_N}.

\bibitem{Madhayama}
U.~K.~V. Sarma, Venketeswara Pai, Dinesh~Mohan Joshi, and K.~Ramasubramanian.
\newblock Madhyam\=anayanaprak\=ara\d h: {A} hitherto unknown manuscript
  ascribed to {M}\=adhava.
\newblock {\em Indian Journal of history of Science}, 46(1):Supplement 1 -- 15,
  2011.
\newblock Available at
  \url{https://cahc.jainuniversity.ac.in/assets/ijhs/Vol46_1_15_Supplement.pdf}.

\bibitem{Catalogue}
K.~Mahadeva Sastri, editor.
\newblock {\em A descriptive catalogue of {S}anskrit manuscripts in the
  {C}urator's {O}ffice {L}ibrary, {T}rivandrum}, volume~IV.
\newblock Govt of H. H. The Maharaja of Travancore, 1939.
\newblock Available at
  \url{https://archive.org/details/in.ernet.dli.2015.490212/}.

\bibitem{Aryabhatiya1930}
K.~Sambasiva Sastri, editor.
\newblock {\em The {\=A}ryabha\d t\=\i ya of {\=A}ryabha\d t\=ac\=arya with the
  {B}h\=a\d sya of {N}\=\i laka\d n\d tha {S}omasutv\=an {P}art {I}. {G}a\d
  nitap\=ada}.
\newblock Trivandrum Sanskrit Series No. 101. Government of Travancore,
  Trivandrum, 1930.
\newblock Available at
  \url{https://archive.org/details/TSS101AryabhatiyaWithTheCommentaryOfNilakantaSomasutvanPart1KSSastri1930_201809}.

\bibitem{YuktiSayahna}
Jye\d s\d thad\=evan.
\newblock {\em Yuktibh\=a\d s\=a {P}art I: {M}athematics (edited with notes by
  Ramavarma (Maru) Thampuran and A. R. Akhileswara Aiyer)}.
\newblock Magalodayam Ltd., Thrissur, Kerala, 1948.
\newblock (This is the first part of the Malayalam book {\em Yuktibh\=a\d s\=a}
  authored by Jye\d s\d thad\=eva.) The book has been reissued as an e-book in
  the pdf format by Sayahna Foundation, Thiruvananthapuram in 2020. The e-book
  is available at \url{https://books.sayahna.org/ml/pdf/yukthibhasha.pdf}.

\bibitem{Sundareswaran2019}
N.~K. Sundareswaran.
\newblock {Nila School of Mathematics - A Note}.
\newblock {\em Calicut University Sanskrit Series}, (59):156--165, 2019.
\newblock (A paper presented in the National Seminar on Regional Traditions of
  Sanskrit: Contributions of North Kerala organised by the Department of
  Sanskrit, University of Calicut from 27.2.2018 to 1.3.2018.) Available at
  \url{https://www.researchgate.net/publication/377437682_Nila_School_of_Mathematics_-A_Note}.

\bibitem{Whish1834}
Charles~M. Whish.
\newblock \uppercase{XXXIII}. \uppercase{O}n the {H}indu quadrature of the
  circle, and the infinite series of the proportion of the circumference to the
  diameter exhibited in the four \'{S}\=astras, the {T}antra {S}angraham,
  {Y}ucti {B}h\=ashá, {C}arana {P}adhati, and {S}adratnam\=la.
\newblock {\em Transactions of the Royal Asiatic Society of Great Britain and
  Ireland}, 3(3):509--523, 1834.
\newblock Available at \url{https://archive.org/details/jstor-25581775}.

\bibitem{Lagnaprakarana03}
Nagakiran Yelluru and Aditya Kolachana.
\newblock {Geometry of {\em pr\=a\d nakal\=ntara} in the {\em Lagnaprakara\d
  na}}.
\newblock {\em Indian Journal of History of Science}, 58:171--180, Sep. 2023.

\end{thebibliography}
\end{document}